\documentclass{article}

\usepackage{amsmath}
\usepackage{amscd}
\usepackage{amssymb}
\usepackage{mathrsfs}
\usepackage{stmaryrd}
\usepackage{pst-node}
\usepackage{tikz-cd}

\newtheorem{theorem}{Theorem}[section]

\newtheorem{proposition}{Proposition}[section]
\newtheorem{lemma}{Lemma}[section]
\newtheorem{corollary}{Corollary}[section]

\newtheorem{conjecture}{Conjecture}[section]

\newcommand{\del}{\nabla}

\def\be{\begin{equation}}
\def\ee{\end{equation}}

\def\endproof{{\hfill $\square$}\medskip}

\begin{document}
\title{Symplectic aspects of polar actions}
\author{XIAOYANG CHEN\footnotemark[1] AND JIANYU OU\footnotemark[2]}
\renewcommand{\thefootnote}{\fnsymbol{footnote}}
\footnotetext[1]{School of Mathematical Sciences, Institute for Advanced Study, Tongji University, Shanghai, China. email: $xychen100@tongji.edu.cn$. %Xiaoyang Chen is partially supported by Scientific Research Foundation 172119, Institute for Advanced Study, Tongji University.
}

\footnotetext[2]{Department of Mathematics, University of Macau, China. email: $eyes\_loki@hotmail.com$. %Jianyu Ou is partially supported by the Project MYRG2015-00235-FST.
}

\maketitle

\begin{abstract}
An isometric compact group action $G \times (M,g)\rightarrow (M,g)$ is called polar if there exists a closed embedded submanifold $\Sigma \subseteq M$ which meets all the orbits orthogonally. Let $\Pi$ be the associated generalized Weyl group. We study the properties of the lifting action $G$ on the cotangent bundle $T^*M$. In particular, we show that the restriction map $(C^\infty(T^*M))^G\rightarrow (C^\infty(T^*\Sigma))^\Pi$ is a surjective homomorphism of Poisson algebras. As a corollary, the singular symplectic reductions ${T^*M}\sslash G$ and $T^*\Sigma\sslash\Pi$ are isomorphic as stratified symplectic spaces, which gives a partial answer to a conjecture of Lerman, Montgomery and Sjamaar.
\end{abstract}

\section{Introduction}

Let $(M,g)$ be a complete Riemannian manifold and $G$ a compact Lie group acting on $M$ by isometries. This action is called polar if there exists a closed embedded submanifold $\Sigma \subseteq M$ meeting all orbits orthogonally (\cite{PalaisTerng}). Then $M$ is called a polar G-manifold and such a submanifold $\Sigma$ is called a section and comes with a natural action by a discrete group of isometries $\Pi = \Pi(\Sigma)$, called its generalized Weyl group. Recall that by definition, $\Pi(\Sigma):= N(\Sigma)/Z(\Sigma)$, where

\begin{eqnarray*}
N(\Sigma)&=&\{g\in G| g\Sigma=\Sigma \},\\
Z(\Sigma)&=&\{g\in G| gx=x,\ x\in\Sigma \}.
\end{eqnarray*}

Polar actions have nice properties and have been studied by many people, see for instance \cite{Dadok}, \cite{GZ}, {\cite{M2}}, {\cite{PalaisTerng}}, \cite{PT}.
Basic examples of polar actions are the adjoint action of a compact Lie group on its Lie algebra. More generally, isotropy representations of symmetric spaces are also polar. It's a classical theorem of Dadok \cite{Dadok} which shows that a polar representation is (up to orbit equivalence) the isotropy representation of a symmetric space. An important feature of polar actions is the following Chevalley Restriction Theorem \cite{PalaisTerng}.

\begin{theorem}
Let $(M,g)$ be a polar G-manifold with a section $\Sigma$ and generalized Weyl group $\Pi$. Then the following restriction to $\Sigma$ is an isomorphism:
$$|_\Sigma:\ \ C^\infty(M)^G \rightarrow C^\infty(\Sigma)^\Pi,$$
where $C^\infty(M)^G$ is the algebra of G-invariant smooth functions on $M$.

\end{theorem}

For a generalization of Chevalley Restriction Theorem to tensors, see \cite{M2}.

In this paper we study symplectic aspects of polar actions. More precisely, we are looking at the lifting action of $G$ on the cotangent bundle $T^*M$ with its canonical symplectic structure $\omega$. This action is a Hamiltonian action with a moment map given by $u: T^*M\rightarrow \mathfrak{g}^*$ with
\begin{eqnarray}
u_X(x,\xi) = \langle \xi, X^*(x)\rangle, {\label{1.1}}
\end{eqnarray}
where $\mathfrak{g}$ is the Lie algebra of $G$, $\mathfrak{g}^*$ is the dual of $\mathfrak{g}$ and $u_X(x,\xi)=\langle u$, $X\rangle(x,\xi),\ X\in \mathfrak{g}$. Moreover, $X^*$ is the vector field on $M$ generated by $X$.
The moment map satisfies the following equations:
\begin{eqnarray*}
 du_X&=&i_{X^{\#}}\omega,\\
u(g\cdot (x,\xi))&=&Ad^*_g\cdot u(x,\xi),\ \forall g\in G,
\end{eqnarray*}
where $X^\#$ is the vector field on $T^*M$ generated by $X$.

Our starting point is the following observation.

\begin{proposition}{\label{prop1}}
Let $(M,g)$ be a polar $G$-manifold with a section $\Sigma$. Then $T^*\Sigma$ meets all $G$-orbits of the action $G\times u^{-1}(0)\rightarrow u^{-1}(0)$.
\end{proposition}

Here $T^*\Sigma$ is seen as a submanifold of $T^*M$ under the natural isomorphism $T^*M \cong TM$ induced by the Riemannian metric $g$.
Note that in general $T^*\Sigma$ can not meet all orbits of the lifting action of $G$ on $T^*M$ as it is easy to see that $T^*\Sigma \subseteq u^{-1}(0)$ from ({\ref{1.1}}).

Let $C^\infty(T^*M)^G$ be the algebra of $G$-invariant smooth functions on $T^*M$. It carries a natural Poisson algebra structure with Poisson bracket $\{f,g\}:=\omega(X_f,X_g)$, where $X_f$ is the Hamiltonion vector field of $f$ satisfying $i_{X_f}\omega=df$.
The generalized Weyl group $\Pi$ of the $G$-action on $M$ also acts on $T^*\Sigma$.
Let $C^\infty(T^*\Sigma)^\Pi$ be the algebra of $\Pi$-invariant smooth functions on $T^*\Sigma$ with its natural Poisson algebra structure. Our main result is the following symplectic analogue of Chevalley Restriction Theorem.

\begin{theorem}{\label{mainthm}}
Let $(M,g)$ be a polar $G$-manifold with a section $\Sigma$ and generalized Weyl group $\Pi$. Then the following restriction to $T^*\Sigma$ is a surjective homomorphism of Poisson algebras:
$$|_{T^*\Sigma}:\ C^\infty(T^*M)^G\rightarrow C^\infty(T^*\Sigma)^\Pi.$$
\end{theorem}

The above restriction to $T^*\Sigma$ is not injective in general as $T^*\Sigma \subseteq u^{-1}(0) \neq T^*M$ unless $u\equiv 0$.

The symplectic reduction $T^*M\sslash G:=u^{-1}(0)/G$ is not a smooth manifold in general. However, it's a stratified symplectic space defined in \cite{SL}. The reader is referred to \cite{SL} for the precise definition of stratified symplectic spaces.
A basic example is given by 
$$X_0=X\sslash G:= J^{-1}(0)/G,$$
where $X$ is a Hamiltonian $G$-space with a moment map $J: X\rightarrow \mathfrak{g}^*$. Following \cite{SL}, we define a function $f_0: X_0 \rightarrow \mathbb{R}$ to be smooth if there exists a function $F\in C^\infty(X)^G$ with $F|_{J^{-1}(0)}=\pi^*f_0$, where $\pi: J^{-1}(0)\rightarrow J^{-1}(0)/G$ is the projection map.
In other words, $C^\infty(X_0)$ is isomorphic to $C^\infty(X)^G/{I^G}$, where $I^G$ is the ideal of $G$-invariant smooth functions on $X$ vanishing on $J^{-1}(0)$. The algebra $C^\infty(X_0)$ inherits a Poisson algebra structure from $C^\infty(X)$.

Let $G$ and $H$ be Lie groups and $X$, resp. $Y$, be smooth manifolds on which $G$, resp. $H$ act properly. The stratified symplectic spaces $T^*X\sslash G$ and $T^*Y\sslash H$ are isomorphic if there exists a homeomorphism $\phi: T^*X\sslash G \rightarrow T^*Y\sslash H$
and the pullback map

\begin{eqnarray*}
\phi^*:\ C^\infty(T^*Y\sslash H) &\rightarrow& C^\infty(T^*X\sslash G),\\
f &\mapsto& f\circ \phi
\end{eqnarray*}
is an isomorphism of Poisson algebras.

In \cite{LMS} (Page 13, Conjecture 3.7), they made the following conjecture.

\begin{conjecture}{\label{conj}}
Let $G$ and $H$ be Lie groups and $X$, resp. $Y$ be smooth manifolds on which $G$, resp. $H$ act properly. Assume that the orbit spaces $X/G$ and $Y/H$ are diffeomorphic in the sense that there exists a homeomorphism $\phi: X/G\rightarrow Y/H$ such that the pullback map $\phi^*$ is an isomorphism from
$C^\infty(Y/H):= C^\infty(Y)^H$ to $C^\infty(X/G):= C^\infty(X)^G$. Then $T^*X\sslash G$ and $T^*Y\sslash H$ are isomorphic.
\end{conjecture}

Using Theorem $\ref{mainthm}$, we can give a partial answer to 
conjecture $\ref{conj}$. More precisely, we have the following corollary.

\begin{corollary}{\label{coro1}}
Let $M$ be a polar $G$-manifold with a section $\Sigma$ and generalized Weyl group $\Pi$. Then $T^*M\sslash G$ and $T^*\Sigma\sslash \Pi$ are isomorphic.
\end{corollary}

Under a slightly different assumption, it was proved that $T^*M\sslash G$ is homomorphic to $T^*\Sigma\slash \Pi$ in \cite{LMS} (Proposition 3.8).

Proposition $\ref{prop1}$, Theorem $\ref{mainthm}$ and Corollary $\ref{coro1}$ will be proved in section 3.
A main ingredient of the proof is a characterization of symplectic slice representations of the lifting action $G$ on $T^*M$, which is done by using the natural Sasaki metric on $T^*M$. Then combining the Multi-variable Chevalley restriction theorem proved by Tevelev \cite{Tev} and other things, we are able to prove our results. For details, see section 3.

\section*{Acknowledgements} The first author is partially supported by Scientific Research Foundation 172119, Institute for Advanced Study, Tongji University. The second author is partially supported by the Project MYRG2015-00235-FST of the University of Macau. Part of this work was done when both authors were visiting the Institute of Mathematical Sciences in the Chinese University of Hong Kong. We also thanks Professors Huai-Dong Cao and Naichung Conan Leung for helpful discussion.

%%%%%%%%%%%%%%%%%%%%%%%%%%%%%%%%%%%%%%%%%%%%%%%%%%%%%%%%%%%%%%%%
\section{Sasaki metrics on $TM$ and $T^*M$.}

In this section we describe the Sasaki metrics on $TM$ and $T^*M$ constructed in {\cite{SA}}. Given a Riemannian metric $g$ on $M$, its Levi-Civita connection determines a splitting $TTM= \mathcal{H}M \oplus \mathcal{V}M$, where $\mathcal{V}M= \ker d\pi$, $\pi: TM \rightarrow M$ is the projection and $\mathcal{H}M$ is spanned by $X^h$, $X$ is a smooth vector field on $M$.
To describe $X^h$, let $(x,v)\in TM$ and $\gamma(t): [0,1] \rightarrow M$ be a smooth curve such that $\gamma(0)=x$, $\gamma'(0)=X(x)$. Let $Y(t) \in T_{\gamma(t)}M$ such that
\begin{equation*}
\begin{cases}
\del_{\gamma'}Y=0\\
Y(0)=v.
\end{cases}
\end{equation*}
Then $X^h(x,v)=:\bar{\gamma}'(0)$, where $\bar{\gamma}\left(t\right)=\left(\gamma\left(t\right), Y\left(t\right)\right)$.
From the definition of $X^h$, we see that $d\pi (X^h(x,v))=X(x)$.
Let $I_g$ be the natural isomorphism $T^*M\cong TM$ induced by the Riemannian metric $g$.
Then using the splitting $TTM= \mathcal{H}M \oplus \mathcal{V}M{\cong}TM\oplus TM$, we define the Sasaki metric $\tilde{g}$ by
$$\tilde{g}\langle\left(X_1,X_2\right), \left(Y_1,Y_2\right)\rangle:=g\langle X_1,Y_1\rangle+g\langle X_2,Y_2\rangle.$$

Define an almost complex structure $J$ by setting $J(X, Y)=(-Y,X).$ Then $\tilde{g}\langle J\cdot, J\cdot\rangle=\tilde{g}\langle\cdot,\cdot\rangle$ and the symplectic form $\Omega := \tilde{g}(J\cdot,\cdot)$ is nothing but the pullback of $\omega$ by the isomorphism $I^{-1}_g: TM\cong T^*M$, where $\omega$ is the standard symplectic form on $T^*M$.

The Sasaki metric on $T^*M$ is the pullback of $\tilde{g}$ under the isomorphism $I_g: T^*M\rightarrow TM$.
The following Lemma will be important for us.

\begin{lemma}{\label{lemma2.1}}
If $\Sigma$ is a totally geodesic submanifold of $(M,g)$, then $T\Sigma$ is a totally geodesic submanifold of $(TM, \tilde{g})$, where $\tilde{g}$ is the Sasaki metric on $TM$.
\end{lemma}

$\bf{Proof:}$ Let $X$ be a smooth vector field on $M$ such that $X(x)\in T_x\Sigma$, $\forall x\in \Sigma$. As $\Sigma$ is totally geodesic, we see that $X^h|_{T\Sigma}$ is a smooth vector field on $T\Sigma$ from the construction of $X^h$.

The vector field $X$ also induces a vertical vector field $X^{\perp}$ on $TM$. We choose a local coordinate to describe $X^{\perp}$. Let $(x^1,...,x^n)$ be a local coordinate system at $x\in M$, where $n=\dim M$. Then any tangent vector $v\in T_xM$ can be decomposed as $v=v^i{\partial\over\partial x_i}$.
The set of parameters $\{x^1,...,x^n,v^1,...,v^n\}$ forms a natural coordinate system of $TM$. The natural frame in $T_{(x,v)}TM$ is given by $\tilde{\partial}_i={\partial\over\partial{x_i}}$ and $\tilde{\partial}_{n+i}={\partial\over\partial{v_i}}$.
Now if $X= X^i{\partial\over{\partial x_i}}$ is a vector field on $M$, then the vertical vector field $X^{\perp}$ on $TM$ is given by $X^{\perp}=X^i\tilde{\partial}_{n+i}$. As $X(x)\in T_x\Sigma,\ \forall x\in \Sigma$, by definition we see that $X^{\perp}|_{T\Sigma}$ is a vector field on $T\Sigma$.

To see that $T\Sigma$ is totally geodesic in $TM$, choose two vector fields $X$, $Y$ on $M$ such that $X(x),\ Y(x)\in T_x\Sigma,\ \forall x\in \Sigma$, then we have the following formula {\cite{Ko}}:

\begin{eqnarray} 
\tilde{\del}_{X^{\perp}} Y^{\perp}&=&0,{\label{sasaki1}}\\
(\tilde{\del}_{X^h}Y^{\perp})(x,v)&=& (\del_X Y)^{\perp}(x,v)+{1\over2}R_x(v, Y_x, X_x)^h(x,v),{\label{sasaki2}}\\
(\tilde{\del}_{X^{\perp}}Y^h)(x,v)&=& {1\over2} (R_x(v,X_x,Y_x))^h(x,v),{\label{sasaki3}}\\
(\tilde{\del}_{X^h}Y^h)(x,v)&=& (\del_X Y)^h(x,v)-{1\over2}R_x(X_x,Y_x,v)^{\perp}(x,v),{\label{sasaki4}}
\end{eqnarray}
where $(x,v)\in T\Sigma$ and $\del$, resp. $\tilde{\del}$ are Levi-Civita connections of $g$, resp. $\tilde{g}$ and $R$ is the Riemann curvature tensor of $g$.

Since $\Sigma$ is totally geodesic, then 
$\del_X Y(x),\ R_x(v,X_x,Y_x) \in T_x\Sigma.$
From ({\ref{sasaki1}}) $-$ ({\ref{sasaki4}}), it follows that $T\Sigma$ is totally geodesic.

\endproof

%%%%%%%%%%%%%%%%%%%%%%%%%%%%%%%%%%%%%%%%%%%%%%%%
\section{A symplectic analogue of Chevalley Restriction Theorem}

In this section we prove Propositon {\ref{prop1}}, Theorem {\ref{mainthm}} and Corollary {\ref{coro1}}. A crucial property of polar actions we will use in the proof is the following result ({\cite{PalaisTerng}} Theorem 4.6):

\begin{proposition}{\label{prop3}}
 Let $M$ be a polar $G$-manifold with a section $\Sigma$. Then the slice representation at $x$ is polar with a section $T_x \Sigma$, $\forall x \in \Sigma$.
\end{proposition}

Given Propositon {\ref{prop3}}, we can now give a proof of Propositon {\ref{prop1}}.

Recall that $u: T^*M\rightarrow \mathfrak{g}^*$ is given by
$$u_X(x, \xi)= \langle\xi, X^*(x)\rangle.$$
Then for any $(x,\xi)\in u^{-1}(0)$, $\langle \xi, X^*(x)\rangle=0$, $\forall X \in \mathfrak{g}$.
Under the isomorphism $I_g: T^*M \cong TM$ induced by the Riemannian metric $g$, the the vector $\xi^\#:= I_g(\xi) \perp T_x(G\cdot x)$, i.e. $\xi^\# \in T_x(G\cdot x)^{\perp}$.

As the isometric action $G\times M \rightarrow M$ is polar with a section $\Sigma$, there exists $ h_1 \in G$ such that $h_1 x\in \Sigma$.
Then $h_1\xi^\# \in T_{h_1 x }(G\cdot x)^{\perp}$.

By Proposition $\ref{prop3}$, the slice representation
$$G_{h_1 x} \times T_{h_1 x }(G\cdot x)^{\perp} \rightarrow T_{h_1 x }(G\cdot x)^{\perp}$$
is polar. Hence there exists $h_2\in G_{h_1 x}$ such that 
$h_2(h_1\xi^\#)\in T_{h_1 x}\Sigma$.

Let $h=h_2 h_1$, then
$$h(x,\xi^\#)=(h x, h{\xi^\#})=(h_1 x, h_2 h_1 {\xi^\#})\in T\Sigma.$$
So $T^*\Sigma$ meets all orbits of the action $G\times u^{-1}(0)\rightarrow u^{-1}(0).$

We proceed to give a proof of Theorem {\ref{mainthm}}. Recall that we have a splitting $TT^*M\cong \mathcal{H}M\oplus\mathcal{V}M$, which induces an isomorphism $TT^*M \overset{(d\pi, I_g)}{\cong} TM \oplus TM$, where $d\pi$ is the differential of the projection $T^*M\rightarrow M$ and $I_g$ is the natural isomorphism $T^*M\cong TM$ induced by the Riemannian metric $g$.

Let $\{x^1,...,x^n, \xi_1,...\xi_n\}$ be a local coordinate of $T^*M$ at $(x,\xi)$ and $\Gamma^k_{ij}$ be the Christoffel symbols of the Levi-Civita connection $\del$ induced by $g$.
Then the horizontal lift of $\partial\over{\partial x_i}$ at $(x,\xi)$ is given by 
$${\tilde{\partial}\over{\partial x_i}}={\partial\over{\partial x_i}}+\Gamma^k_{il}\xi_k{\partial\over{\partial \xi^l}}.$$
Here a horizontal lift of a vector $X$ at $(x,\xi)$ is defined to be the unique vector $\tilde{X}\in \mathcal{H}M$ such that 
$d\pi(x,\xi)(\tilde{X})=X$.

In terms of local coordinate system $\{x^1,...,x^n, \xi_1,...\xi_n\}$, the almost complex structure $J$ defined in section 2 can be rephrased as

\begin{eqnarray*}
J({\tilde{\partial}\over{\partial x_i}})&=&g_{ij}{\partial\over{\partial\xi^j}},\\
J({\partial\over{\partial \xi^i}})&=&-g^{ij}{\tilde{\partial}\over{\partial x_j}},
\end{eqnarray*}
where $g_{ij}=g({\partial\over{\partial x_i}},{\partial\over{\partial x_j}})$ and $(g^{ij})$ is the inverse matrix of $(g_{ij})$.

Let $X^*=X^i{\partial\over{\partial x_i}}$ be a vector field on $M$ generated by $X\in \mathfrak{g}$. Then the corresponding vector field on $T^*M$ generated by $X$ is 
\begin{eqnarray*}
X^\#(x,\xi)=X^i{\partial\over{\partial x_i}}-\sum_{i,j}{\partial X^j\over{\partial x_i}}\xi_j{\partial\over{\partial \xi^i}},
\end{eqnarray*}
see {\cite{BV}} (Page 16 Lemma 11).

The Sasaki metric $\tilde{g}$ on $T^*M$ satisfies
\begin{eqnarray*}
\tilde{g}\langle{\tilde{\partial}\over\partial x_i},{\tilde{\partial}\over\partial x_j}\rangle&=&g_{ij},\\
\tilde{g}\langle{\tilde{\partial}\over\partial x_i},{{\partial}\over\partial \xi^j}\rangle&=&0,\\
\tilde{g}\langle{{\partial}\over\partial \xi^i},{{\partial}\over\partial \xi^j}\rangle&=& g^{ij}.
\end{eqnarray*}

\begin{lemma}
$\forall (x,\xi)\in T^*\Sigma$, the Sasaki metric $\tilde{g}$ on $T^*M$ induces an orthogonal splitting 
$$T_{(x,\xi)}T^*M=T_{(x,\xi)}(G\cdot(x,\xi))\oplus JT_{(x,\xi)}(G\cdot(x,\xi))\oplus V$$
with $T_{(x,\xi)}T^*\Sigma\subseteq V$ and $V$ is the orthogonal complement of $T_{(x,\xi)}(G\cdot(x,\xi))\oplus JT_{(x,\xi)}(G\cdot(x,\xi))$.
\end{lemma}

{\bf{Proof:}} Let $X^\#_i$ be two vector fields on $T^*M$ generated by $X_i\in \mathfrak{g}$, $i=1,2$ respectively and $Y \in T_{(x,\xi)}T^*\Sigma$. Then $\tilde{g}\langle J X^\#_1, X^\#_2\rangle=\omega(X^\#_1,X_2^\#)=(i_{X_1^\#}\omega)(X^\#_2).$
Let $u$ be the moment map defined in ($\ref{1.1}$),
as $(x,\xi)\in T^*\Sigma\subseteq u^{-1}(0)$, by the $G$-equivalence of $u$, we get $G(x,\xi)\subseteq u^{-1}(0)$.
Hence $$\tilde{g}\langle JX_1^\#,X_2^\#\rangle=(i_{X_1^\#}\omega)(X^\#_2)=du_{X_1}(X_2^\#)=0.$$ 

By the definition of $J$, we get $JT_{(x,\xi)}T^*\Sigma\subseteq T_{(x,\xi)}T^*\Sigma.$
As $T^*\Sigma\subseteq u^{-1}(0),$ we get
$$\tilde{g}\langle X^\#_1, Y\rangle=\tilde{g}\langle JX_1^\#,JY\rangle=\omega(X_1^\#, JY)=(i_{X_1^\#}\omega)(JY)=du_{X_1}(JY)=0.$$
Similarly, $\tilde{g}\langle JX_1^\#, Y\rangle=0$. Hence $T_{(x,\xi)}T^*\Sigma \subseteq V.$
\endproof

The representation
$$G_{(x,\xi)} \times V\rightarrow V$$
is called the symplectic slice representation at $(x,\xi)$.
Note that $G_{(x,\xi)}=(G_x)_\xi=:\{h\in G_x|\ h\xi=\xi\}$.

The following Lemma will be crucial for us.

\begin{lemma}{\label{lemma3.1}}
Let $M$ be a polar $G$-manifold with a section $\Sigma$. Then the symplectic slice representation at $(x,\xi)\in T^*\Sigma$ is the diagonal action (up to identification)
$$(G_x)_{\xi^\#}\times (W\oplus W)\rightarrow W\oplus W,$$
where $W := (G_x \xi^\#)^{\perp}$ is the orthogonal complement of $G_x \xi^\#$ in the slice $(G\cdot x)^{\perp},$ i.e. we have
\begin{eqnarray*}
T_xM&=&T_x(G\cdot x)\oplus (T_x(G\cdot x)^{\perp}),\\
T_x(G\cdot x)^{\perp}&=& G_x \xi^\# \oplus (G_x \xi^\#)^{\perp}.
\end{eqnarray*}
\end{lemma}

${\bf{Proof:}}$ Let $G_{(x,\xi)}\times V\rightarrow V$ be the symplectic slice representation at $(x,\xi)$. Under the isomorphism $\Phi: \mathcal{H}M\oplus \mathcal{V}M \overset{(d\pi,I_g)}{\longrightarrow}TM\oplus TM$, we first claim that
$$\Phi(V)=W\oplus W.$$

Choose a local coodinate system $\{x^1,...,x^n,\xi_1...,\xi_n\}$
of $T^*M$ at $(x,\xi)$. Then we have
\begin{eqnarray*}
d\pi({\tilde{\partial}\over\partial x_i})&=&{\partial\over \partial x_i},\\
I_g({\partial\over\partial\xi^i})&=& g^{ij}{\partial\over \partial x_j}.
\end{eqnarray*}
Let $Z= a^i{\tilde{\partial}\over\partial x_i}+b_i{\partial\over\partial\xi^i}\ \in TT^*M$. Then $\Phi(Z)=(d\pi,I_g)(Z)=(a^i{\partial\over\partial x_i},g^{ij}b_i{\partial\over\partial x_j})=:(Y_1,Y_2).$

Let $X^*=X^i{\partial\over\partial x_i}$ be the vector field on $M$ generated by $X\in\mathfrak{g}$, then the corresponding vector field on $T^*M$ is 
\begin{eqnarray*}
X^\#(x,\xi)&=& X^i{\partial\over\partial x_i}-\sum_{i,j}{\partial X^j\over\partial x_i}\xi_j{\partial\over\partial \xi^i}\\
&=&X^i{\tilde{\partial}\over\partial x_i}-X^i\Gamma^k_{il}\xi_k{\partial\over\partial\xi^l}-\sum_{ij}{\partial X^j\over\partial x_i}\xi_j{\partial\over\partial \xi^i}\\
&=&X^i{\tilde{\partial}\over\partial x_i}-g\langle{\del_{\partial\over\partial x_i}}X^*,\xi^\#\rangle{\partial \over \partial \xi^i}.
\end{eqnarray*}

Then we have

\begin{eqnarray}
\tilde{g}(X^\#,Z)&=&\tilde{g}\langle X^i{\tilde{\partial}\over\partial x_i}-g\langle{\del_{\partial\over\partial x_i}}X^*,\xi^\#\rangle{\partial \over \partial \xi^i}, a^j{\tilde{\partial}\over\partial x_j}+b_j{\partial\over\partial\xi^j}\rangle{\nonumber}\\
&=& g_{ij}X^ia^j-g^{ij}b_j g\langle{\del_{\partial\over\partial x_i}}X^*,\xi^\#\rangle{\nonumber}\\
&=&g\langle X^*, Y_1\rangle- g\langle \del_{Y_2}X^*, \xi^\#\rangle. {\label{3.1}}
\end{eqnarray}

We also have

\begin{eqnarray}
\tilde{g}\langle X^\#, JZ\rangle&=& \tilde{g}\langle X^i{\tilde{\partial}\over\partial x_i}-g\langle{\del_{\partial\over\partial x_i}}X^*,\xi^\#\rangle{\partial \over \partial \xi^i}, J(a^j{\tilde{\partial}\over\partial x_j}+b_j{\partial\over\partial\xi^j})\rangle{\nonumber}\\
&=&\tilde{g}\langle X^i{\tilde{\partial}\over\partial x_i}-g\langle{\del_{\partial\over\partial x_i}}X^*,\xi^\#\rangle{\partial \over \partial \xi^i}, a^jg_{jk}{{\partial}\over\partial \xi^k}-b_j g^{jk}{\tilde{\partial}\over\partial x_k}\rangle{\nonumber}\\
&=&-X^ib_i-a^i g\langle\del_{\partial\over \partial x_i}X^*,\xi^\#\rangle{\nonumber}\\
&=&-g\langle X^*, Y_2\rangle- g\langle\del_{Y_1}X^*,\xi^\#\rangle.
\end{eqnarray}

Now we proceed to prove $\Phi(V)=W\oplus W$.
Let $Z\in TT^*M$ such that $\Phi(Z)=(Y_1,Y_2)\in W\oplus W$.
We claim that $Z\in V$ and it follows that $W\oplus W\subseteq \Phi(V).$
In fact, as $(Y_1,Y_2)\in W\oplus W$, we get
\begin{eqnarray}
g\langle X^*, Y_1\rangle=0, {\label{3.3}}\\
g\langle X^*, Y_2\rangle=0. {\label{3.4}}
\end{eqnarray}

As $X^*$ is a Killing vector field, we get
\begin{eqnarray}
g\langle \del_{Y_2}X^*, \xi^\#\rangle=-g\langle\del_{\xi^\#}X^*, Y_2\rangle=g\langle\del_{\xi^\#}Y_2, X^*\rangle, {\label{3.5}}
\end{eqnarray}
and 
\begin{eqnarray}
g\langle \del_{Y_1}X^*, \xi^\#\rangle=g\langle \del_{\xi^\#}Y_1,X^*\rangle. {\label{3.6}}
\end{eqnarray}

As $M$ is a polar $G$-manifold with a section $\Sigma$, By Proposition {\ref{prop3}}, the slice representation $G_x \times T_x(G\cdot x)^{\perp}\rightarrow T_x(G\cdot x)^{\perp}$ is polar with a section $T_x\Sigma$. Then by Proposition {\ref{prop3}} again, the slice representation
$(G_x)_{\xi^\#} \times W \rightarrow W$ is polar with a section $T_{\xi^\#}T_x\Sigma$.
As $Y_1\in W$, there exists $h\in (G_x)_{\xi^\#}$ such that $hY_1\in T_{\xi^\#}T_x\Sigma$. Hence $Y_1\in h^{-1}(T_{\xi^\#}(T_x\Sigma))=T_{\xi^\#}T_x(h^{-1}\Sigma)\cong T_x(h^{-1}\Sigma).$ We also have $\xi^\#=h^{-1}\xi^\#\in T_x(h^{-1}\Sigma),$ as $\Sigma$ is totally geodesic ({\cite{PalaisTerng}}, Theorem 3.2), so is $h^{-1}\Sigma$.

Then
\begin{eqnarray}
g\langle \del_{\xi^\#}Y_1, X^*\rangle=g\langle B(\xi^\#,Y_1),X^*\rangle=0 {\label{3.7}},
\end{eqnarray}

\begin{eqnarray}
g\langle \del_{\xi^\#}Y_2,X^*\rangle = g\langle B(\xi^\#, Y_2), X^*)=0 {\label{3.8}}
\end{eqnarray}
where $B(\ ,\ )$ is the second fundamental form of $h^{-1}\Sigma$.

By $(\ref{3.1}),\ (\ref{3.3}),\ (\ref{3.5})$ and $(\ref{3.8})$, we get
$$\tilde{g}\langle X^\#, Z\rangle=g\langle X^*, Y_1\rangle-g\langle\del_{Y_2}X^*,\xi^\#\rangle=0.$$

Similary we get
$\tilde{g}\langle X^\#, JZ\rangle=0$. Hence $\tilde{g}\langle JX^\#, Z\rangle=-g\langle X^\#, JZ\rangle=0$.
It follows that $Z\in V$, which implies that $W\oplus W\subseteq \Phi(V)$.

On the other hand, we claim that $\dim{(W\oplus W)}= \dim \Phi(V)$.
In fact,
\begin{eqnarray*}
\dim (W\oplus W) &=& 2\dim W\\
&=& 2(\dim (T_x G\cdot x)^{\perp}-\dim(G_x\cdot \xi^\#))\\
&=& 2(\dim M- \dim G\cdot x-(\dim G_x-\dim(G_x)_{\xi^\#}))\\
&=& 2 \dim M- 2(\dim G- \dim G_{(x,\xi)}),
\end{eqnarray*}
and 
\begin{eqnarray*}
\dim \Phi(V)&=& \dim V\\
&=& \dim T^*M- 2 \dim G\cdot(x,\xi)\\
&=& 2(\dim M- \dim G\cdot(x,\xi))\\
&=& 2 (\dim M- (\dim G-\dim G_{(x,\xi)})).
\end{eqnarray*}
Hence $\dim \Phi(V)= \dim{(W\oplus W)}$ and we have $\Phi(V)= W\oplus W$.

Now Lemma {\ref{lemma3.1}} follows from the following commutative diagram

\[\begin{tikzcd}
V \arrow{r}{\Phi} \arrow[swap]{d}{G_{(x,\xi)}} & {W\oplus W} \arrow{d}{(G_x)_{\xi^\#}} \\
V \arrow{r}{\Phi} & W\oplus W
\end{tikzcd}
\]

\endproof

Given Lemma {\ref{lemma3.1}}, we can now give a proof of Theorem {\ref{mainthm}}. We first show that the following restriction map is surjective:
\begin{eqnarray*}
|_{T^*\Sigma}:\ C^\infty(T^*M)^G\rightarrow C^\infty(T^*\Sigma)^\Pi.
\end{eqnarray*}

$\forall (x,\xi)\in T^*\Sigma$, the Sasaki metric $\tilde{g}$ on $T^*M$ induces an orthogonal splitting
$$T_{(x,\xi)}T^*M=T_{(x,\xi)}G(x,\xi)\oplus J T_{(x,\xi)}G(x,\xi)\oplus V,$$
where $\Phi(V)\cong W\oplus W$ by Lemma {\ref{lemma3.1}}.

The Slice Theorem says that for an open $G$-invariant tubular neighborhood $U_{(x,\xi)}$ of the orbit $G(x,\xi)$, there is a $G$-equivalent diffeomorphism
\begin{eqnarray*}
\exp^{\perp}:\  G\times_{G_{(x,\xi)}}S^{\perp}_{(x,\xi)}(\epsilon)\rightarrow U_{(x,\xi)},
\end{eqnarray*}
where $S^{\perp}_{(x,\xi)}:= JT_{(x,\xi)}G(x,\xi)\oplus V$,
$S^{\perp }_{(x,\xi)}(\epsilon)$ is the $\epsilon$-ball in $S^{\perp}_{(x,\xi)}$
 and $\exp^{\perp}$ is the normal exponential map of $G(x,\xi)$.

Let $U= \underset{(x,\xi)\in T^*\Sigma}{\bigcup} U_{(x,\xi)}$. As $T^*\Sigma$ intersects all orbits in $u^{-1}(0)$ by Proposition {\ref{prop1}}, we see that $U$ is a $G$-invariant open neighborhood of $u^{-1}(0)$.
$\forall f\in C^\infty(T^*\Sigma)^\Pi$, we first show that there exists
$F_{\epsilon}\in C^\infty(U_{(x,\xi)})^G$ such that
\begin{eqnarray}
F_{\epsilon}|_{T^*\Sigma \cap U_{(x,\xi)}}=f|_{T^*\Sigma \cap U_{(x,\xi)}} {\label{3.7}}
\end{eqnarray}

By the existence of $G$-invariant partition of unity subject to the cover 
$U= \underset{(x,\xi)\in T^*\Sigma}{\bigcup} U_{(x,\xi)}$,
then there exists $F\in C^\infty(U)^G$ such that $F|_{T^*\Sigma}=f$.
Extending $F$ to $\tilde{F}\in C^\infty(T^*M)^G$, we then prove our desired result.

To prove $(\ref{3.7})$, we first recall some facts on polar representations which we will use. Let $(G, K)$ be a symmetric pair and consider the isotropy representation of $K$ on $\mathfrak{p}= T_K{(G\slash K)}$. It is a polar action and any maximal abelian sub-algebra $\Sigma$ is a section. Its generalized Weyl group $\Pi$ is also called the "baby" Weyl group. Consider the diagonal action of $K$ on $\mathfrak{p}^m$ (respectively $\Pi$ on $\Sigma^m$) and the corresponding algebra of invarant ($m$-variable) polynomials $\mathbb{R}[\mathfrak{p}^m]^K$ (respectivelly $\mathbb{R}[\Sigma^m]^\Pi$). Then we have the following result due to Tevelev {\cite{Tev}}.

\begin{theorem}{\label{thm3.1}}
The restriction map $|_{\Sigma}:\ \mathbb{R}[\mathfrak{p}^m]^K\rightarrow \mathbb{R}[\Sigma^m]^\Pi$ is surjective.
\end{theorem}

As a polar representation is (up to orbit equivalence) the isotropy representation of a symmetric space {\cite{Dadok}}. Theorem {\ref{thm3.1}} generalizes to the class of polar representations {\cite{M2}} (Corollary 2).

\begin{corollary}{\label{coro2}}
Let $K\subseteq O(\mathfrak{p})$ be a linear representation which is also polar with a section $\Sigma$ and generalized Weyl group $\Pi$. Then the restriction is surjective:
$$|_{\Sigma}:\ \mathbb{R}[\mathfrak{p}^m]^K\rightarrow \mathbb{R}[\Sigma^m]^\Pi.$$
\end{corollary}

\begin{corollary}{\label{coro3}}
Let $\mathfrak{p}$ be a polar representation of a compact Lie group $K$ with a section $\Sigma$ and generalized Weyl group $\Pi$. Then the restriction to $\Sigma$ is surjective:
$$|_{\Sigma}:\ C^\infty(\mathfrak{p}^m)^K\rightarrow C^\infty(\Sigma^m)^\Pi.$$
\end{corollary}

{\bf{Proof:}} It's a classical result of Hillbert ({\cite{Sp}}, Proposition 2.4.14) that $\mathbb{R}[\mathfrak{p}^m]^K$ is finitely generated. Let $\rho_1,...,\rho_n$ be generators. By Corollary {\ref{coro2}}, $\rho_1|_{\Sigma},...,\rho_n|_{\Sigma}$ generate $\mathbb{R}[\Sigma^m]^\Pi$.

For any $f\in C^\infty(\Sigma^m)^\Pi$, apply Schwarz's Theorem {\cite{S}} to the action of $\Pi$ on $\Sigma^m$, we get $F\in C^\infty(\mathbb{R}^n)$ such that $f=F\circ \rho|_{\Sigma}$, where
$\rho|_{\Sigma}:\ \Sigma^m \rightarrow \mathbb{R}^n$ be the map whose coordinates are $\rho_1|_{\Sigma},...,\rho_n|_{\Sigma}$.
Then $\tilde{f}=F\circ \rho \in C^\infty(\mathfrak{p}^m)^K$ such that $\tilde{f}|_{\Sigma}=f$.
 
\endproof

We can now give a proof of ({\ref{3.7}}). By Lemma {\ref{lemma2.1}}, as $\Sigma$ is totally geodesic in $M$, then $T^*\Sigma$ is totally geodesic in $T^*M$. Hence the normal exponential map $\exp^{\perp}$ of the orbit $G\cdot(x,\xi)$ maps the $\epsilon$-ball $B_{\epsilon}$ in $T_{(x,\xi)}T^*\Sigma{\cong}T_x\Sigma\oplus T_x\Sigma$ diffeomorphically onto $T^*\Sigma \cap U_{(x,\xi)}$.
$\forall f \in C^{\infty}(T^*\Sigma)^\Pi$, $f\circ \exp^{\perp}:\ B_{\epsilon}\rightarrow \mathbb{R}$ is a $\Pi_{(x,\xi)}$-invariant smooth function, where $\Pi_{(x,\xi)}=\{h\in\Pi|\ h(x,\xi)=(x,\xi)\}$. Let $W$ be a polar representation of $K:= G_{(x,\xi)}$ with a section $T_{\xi^\#}T_x\Sigma\cong T_x\Sigma$ defined in Lemma {\ref{lemma3.1}}. By Corollary {\ref{coro3}}, we see that there exists $f_{\epsilon}\in C^\infty(W\oplus W)^K$ such that
$$f_{\epsilon}|_{B_{\epsilon}}=f\circ\exp^{\perp}$$
Hence $f_{\epsilon}\circ(\exp^{\perp})^{-1}=f$ on $T^*\Sigma\cap U_{(x,\xi)}$. Combined with Lemma {\ref{lemma3.1}} and the Slice theorem, then $f_{\epsilon}$ is pulled back to be a smooth function on $G\times S^{{\perp}}_{(x,\xi)}(\epsilon)$ which descends to $F_{\epsilon}\in C^\infty(U_{(x,\xi)})^G$ such that $F_{\epsilon}=f$ on $T^*\Sigma\cap U_{(x,\xi)}$. We finish the proof of the surjectivity part of Theorem {\ref{mainthm}}.

Let $\omega$ be the standard symplectic form on $T^*M$. We show that the restriction to $T^*\Sigma$ preserves Poisson brackets 
$(C^\infty(T^*M)^G,\{\ ,\ \}_1)$ and $(C^\infty(T^*\Sigma)^\Pi,\{\ ,\ \}_2)$, where $\{\ ,\ \}_{i=1,2}$ are Poisson brackets induced by $\omega$ and $\omega|_{T^*\Sigma}$ respectively.

Let $\mathring{M}\subseteq M$ be the union of principal orbits and $\mathring{\Sigma}=\Sigma\cap \mathring{M}$.
Then $\mathring{\Sigma}$ is open and dense in $\Sigma$ ({\cite{GZ}} Propsition 1.3). It follows that $T^*\mathring{\Sigma}\subseteq T^*\Sigma$ is also open and dense. $\forall (x,\xi)\in T^*\mathring{\Sigma}$, we have the following orthogonal splitting with respect to the Sasaki metric $\tilde{g}$ on $T^*M$:

\begin{eqnarray}{\label{split}}
T_{(x,\xi)}T^*M\cong T_{(x,\xi)}G(x,\xi)\oplus JT_{(x,\xi)}G(x,\xi)\oplus T_{(x,\xi)}T^*\mathring{\Sigma}.
\end{eqnarray}
To see this, as $\mathring{\Sigma}$ consists of principal orbits, the slice representation at $x\in \mathring{\Sigma}$ is trivial. Hence $G_{(x,\xi)}=G_x$, $\forall(x,\xi)\in T^*\mathring{\Sigma}$. By {\cite{PalaisTerng}}, we also have $\dim G\cdot x+\dim\mathring{\Sigma}=\dim M$.
Then the dimension of the vector space on the right hand side of ({\ref{split}}) is equal to 
$$2\dim G(x,\xi)+2\dim\mathring{\Sigma}=2(\dim G\cdot x+\dim\mathring{\Sigma})=2\dim M$$
which finishes the proof of (\ref{split}).

$\forall f \in C^\infty(M)^G$, at $(x,\xi)\in T^*\mathring{\Sigma}$, we can write $$X_f=X+JY+Z,$$ where $X,Y\in T_{(x.\xi)}G(x,\xi)$, $Z\in T_{(x,\xi)}T^*\mathring{\Sigma}$.

Recall that $i_{X_f}\omega=df$, $\omega$ is the standard symplectic form on $T^*M$.
Since $f$ is $G$-invariant, we get $(i_{X_f}\omega)(Y)=df(Y)=0.$
Then $$\omega(X_f,Y)=\tilde{g}(JX_f,Y)=\tilde{g}(JX-Y+JZ,Y)=-\tilde{g}(Y,Y).$$
It follows that $Y=0$.

Now let $f_1,f_2\in C^\infty(T^*M)^G$, at $(x,\xi)\in T^*\mathring{\Sigma},$ we have 
$$X_{f_i}=X_i+Z_i,\ i=1,2$$
where $X_i\in T_{(x,\xi)}G(x,\xi)$, $Z_i\in T_{(x,\xi)}T^*\mathring{\Sigma}$.

Let $\hat{f}_1={f_1}|_{T^*\mathring{\Sigma}}$,
then we claim that $X_{\hat{f}_1}=Z_1$.
In fact $i_{X_{\hat{f}_1}}\omega|_{T^*\Sigma}=d{\hat{f}_1}$. $\forall Y\in T_{(x,\xi)}T^*\mathring{\Sigma}$, we have
\begin{eqnarray*}
\tilde{g}\langle X_{\hat{f}_1}, Y\rangle&=&\omega(X_{\hat{f}_1},JY)\\
&=&d{\hat{f}_1}(JY)\\
&=&df_1(JY)\\
%&=&(i_{X_{f_1}}\omega)(JY)\\
&=&\omega(X_{f_1}, JY)\\
%&=&\tilde{g}\langle X_{f_1}, Y\rangle\\
%&=&\tilde{g}\langle X_1+Z_1, Y\rangle\\
&=&\tilde{g}\langle Z_1, Y\rangle.
\end{eqnarray*}

Then at $(x,\xi)$, 
\begin{eqnarray*}
\{f_1,f_2\}_1&=&\omega(X_{f_1},X_{f_2})\\
%&=&\tilde{g}(JX_{f_1},X_{f_2})\\
&=&\tilde{g}(JX_1+JZ_1,X_2+Z_2)\\
&=&\tilde{g}(JZ_1,Z_2)\\
&=&\omega|_{T^*\Sigma}(Z_1,Z_2)\\
&=&\{\hat{f}_1,\hat{f}_2\}_2.
\end{eqnarray*}
By continuity, $\{f_1,f_2\}_1=\{\hat{f}_1,\hat{f}_2\}_2$ on 
$T^*\Sigma$ everywhere.

{\bf{Proof of Corollary {\ref{coro1}}}:}
As $M$ is a polar $G$-manifold, the inclusion: $\Sigma\slash \Pi \rightarrow M\slash G$ is a homomophism {\cite{PalaisTerng}}. By Chevalley restriction theorem {\cite{PalaisTerng}}, the restriction 
$|_{\Sigma}: C^\infty(M)^G\rightarrow C^\infty(\Sigma)^\Pi$ is an isomorphism, which implies that $M/G$ is diffeomophic to $\Sigma/\Pi$.
We now show that $T^*M\sslash G$ is diffeomophic to $T^*\Sigma\sslash \Pi$. By Theorem {\ref{mainthm}}, $C^\infty(T^*M)^G/I^G$ is isomophic to $C^\infty(T^*\Sigma)^\Pi$ as Poisson algebra, where $I^G$ is the ideal of $G$-invariant smooth functions on $X$ vanishing on $u^{-1}(0)$.

It's enough to show that the inclusion $T^*\Sigma\sslash\Pi\rightarrow T^*M\sslash G$ is a homomorphism. By Proposition {\ref{prop1}}, it sufficies to show
$$G\cdot(x,\xi)\cap T^*\Sigma=\Pi\cdot(x,\xi),\ \forall (x,\xi)\in T^*\Sigma.$$

Clearly $\Pi\cdot(x,\xi)\subseteq G\cdot(x,\xi)\cap T^*\Sigma.$
On the other hand, $\forall h_1(x,\xi)\in G(x,\xi)\cap T^*\Sigma$, we have $h_1 x \in G\cdot x\cap \Sigma$ and $h_1\xi^\#\in T_{h_1x}\Sigma$.
By Corollary 4.9 in {\cite{PalaisTerng}}, we get
$$G\cdot x\cap\Sigma= \Pi\cdot x,\ \forall x\in\Sigma.$$
Hence
\begin{eqnarray}{\label{h1h2}}
h_1x=h_2x,\ h_2\in \Pi.
\end{eqnarray}
Then $(h_2^{-1}h_1)x=x$ and so $h_2^{-1}h_1\in G_x$.
Since $(x,\xi)\in T^*\Sigma$, we get $\xi^\#\in T_x(G\cdot x)^{\perp}.$
By Lemma {\ref{lemma3.1}}, the slice representation: $G_x\times T_x(G\cdot x)^{\perp}\rightarrow T_x(G\cdot x)^{\perp}$ is polar with a section $T_x\Sigma$ and generalized Weyl group $\Pi_x$. By Corollary 4.9 in {\cite{PalaisTerng}} again,
$$G_x\cdot\xi^\#\cap T_x\Sigma=\Pi_x\cdot\xi^\#.$$
As $h_2^{-1}h_1\in G_x$, $h_2\in \Pi$, $h_1\xi^\#\in T_{h_1x}\Sigma$, we get $h_2^{-1}h_1\xi^\#\in G_x\cdot\xi^\#\cap T_x\Sigma.$
Then there exists $h_3\in \Pi_x$ such that
$$h_2^{-1}h_1\xi^\#=h_3\xi^\#.$$
Hence $h_1\xi^\#=h_2h_3\xi^\#\in\Pi\cdot\xi^\#.$
Combined with {\ref{h1h2}}, we obtain
$h_1(x,\xi^\#)=(h_2x,h_2h_3\xi^\#)=h_2h_3(x,\xi^\#)\in \Pi\cdot(x,\xi^\#).$
So $G(x,\xi)\cap T^*\Sigma\subseteq \Pi\cdot(x,\xi).$
\endproof


\begin{thebibliography}{}

%\bibitem{BlumenthalHebda} R.A Blumenthal and J.J. Hebda, \textit{De Rham decomposition theorems for foliated manifolds}, Ann. Inst. Fourier 33 (1983), 183-198.

\bibitem{BV} N. Berline and M. Vergne, \textit{Hamiltonian manifolds and moment map}, http://www.cmls.polytechnique.fr/perso/berline/cours-Fudan.pdf

\bibitem{Dadok} J. Dadok, \textit{Polar coordinates induced by actions of compact Lie groups}, Trans. Amer. Math. Soc. 288(1985), 125-137.

%\bibitem{DK} J. Dadok and V. Kac, \textit{Polar representations}, J. Algebra, 92(1985), 504-524 .

%\bibitem{Delzant} T. Delzant, \textit{Hamiltoniens p$\acute{e}$riodiques et image convexe de l'application moment}, Bull. Soc. Math. France 116 (1988), 315-339.

%\bibitem{De2} T. Delzant, \textit{Classification des actions hamiltoniennes compl$\grave{e}$tement int$\acute{e}$grables de rang deux}, Ann. Global Anal. Geom. 8 (1990), 87-112.

%\bibitem{Do} P. Dombrowski, \textit{On the geometry of the tangent bundle}, J. Reine Angew. Math. 210 (1962), 73-88.


%\bibitem{FG} F. Fang, K. Grove and G. Thorbergsson, \textit{Tits geometry and positive curvature}, to appear in Acta Math.

\bibitem{GZ} K. Grove and W. Ziller, \textit{Polar manifolds and actions}, J. Fixed Point Theory Appl. 11 (2012), no.2, 279-313.


%\bibitem{Guillemin} V. Guillemin, \textit{Moment maps and Combinatorial Invariants of Hamiltonian $T^n$-spaces}, Birkh$\ddot{a}$user, 1994.

%\bibitem{GS} V. Guillemin and S Sternberg, \textit{Multiplicity-free spaces}, J. Diff. Geom. 19 (1984), 31-56.

%\bibitem{Igl} P. Igl$\acute{e}$sias, \textit{Les SO(3)-vari$\acute{e}$t$\acute{e}$s symplectiques et leur classification en dimension 4}, Bull. Soc. Math. France 199(1991), 371-396.


%\bibitem{KarshonTolman01} Y. Karshon and S. Tolman, \textit{Centered complexity one Hamiltonian torus actions.}, Trans. Amer. Math. Soc. 353 (2001), no. 12, 4831-4861.


%\bibitem{KarshonTolman03} Y. Karshon and S. Tolman, \textit{Complete invariants for Hamiltonian torus actions with two dimensional quotients},  J. Symplectic Geom. 2 (2003), no. 1, 25-82.

%\bibitem{KarshonTolman14} Y. Karshon and S. Tolman, \textit{Classification of Hamiltonian torus actions with two-dimensional quotients}, Geom. Topol. 18 (2014), no. 2, 669-716. 

%\bibitem{Kirwan} F. Kirwan, \textit{Cohomology of Quotients in Symplectic and Algebraic Geometry}, Princeton University Press, Princeton, 1984.


%\bibitem{Knop} F. Knop, \textit{Automorphisms of multiplicity free Hamiltonian manifolds}, J. Amer. Math. Soc. 24 (2011), no.2, 567-601.

\bibitem{Ko} O. Kowalski, \textit{Curvature of the induced Riemannian metric on the tangent bundle of a Riemannian Manifold}, J. Reine Angew. Math. 250(1971), 124-129.

%\bibitem{Lerman} E. Lerman, \textit{Symplectic cuts}. Math. Res. Lett. 2 (1995), 247-258.


\bibitem{LMS} E. Lerman, R. Montgomery and R. Sjamaar, \textit{Examples of singular reduction}， in Symplectic geometry, D. A. Salamon, ed., Cambridge: Cambridge University
Press, 1993.


%\bibitem{Los} I. Losev, \textit{Proof of the Knop conjecture}, Ann. Inst. Fourier (Grenoble) 59 (2009), 1105-1134.

%\bibitem{MiFo} A. Miscenko, A. Fomenko, \textit{A generalized Liouville method for the integration of Hamiltonian system}, Funkcional. Anal. i Prilozen. 12 (1978), 46-56, 96.


%\bibitem{M1} R A. E. Mendes, \textit{Equivariant tensors on polar manifolds}, PhD dissertation (2011). 

\bibitem{M2} R. A. E. Mendes, \textit{Extending tensors on polar manifolds}, Math. Ann. 365(2016), no.3, 1409-1424.


\bibitem{PalaisTerng} R. S. Palais and C. L. Terng, \textit{A general theory of canonical forms}, Trans. Amer. Math. Soc. 300 (1987), no. 2, 771-789. 

\bibitem{PT} F. Podest$\grave{a}$, G. Thorbergsson, \textit{Polar actions on rank-one symmetric spaces}, J. Differential Geom. 53 (1999), 131-175.

\bibitem{SA} S. Sasaki, \textit{On the differential geometry of tangent bundles of Riemannian manifolds}. Tohoku Math. J. Vol. 10 (1958), 338-354.

%\bibitem{S0} G.W. Schwarz, \textit{Generalized Orbit Spaces}, revised version of PhD thesis, MIT, 1972, unpublished.

\bibitem{S} G. W. Schwarz, \textit{Smooth functions invariant under the action of a compact Lie group}, Topology, 14 (1975), no 1, 63-68.

\bibitem{SL} S. Sjamaar and E. Lerman, \textit{Stratified symplectic spaces and reduction}, Ann. Math. 134(1991), 375-422.

\bibitem{Sp} T. A. Springer, \textit{Invariant Theory}, Lecture Notes in Mathematics, (1997) Vol. 585. Springer, Berlin.

\bibitem{Tev} E. A. Tevelev, \textit{On the Chevalley restriction theorem.}, J. Lie Theory. 10 (2000), no. 2, 323-330.

%\bibitem{Woo} C. Woodward, \textit{The classification of transversal multiplicity-free group actions}, Ann. Global Anal. Geom. 14 (1996), no. 1, 3-42.

\end{thebibliography}
\end{document}